\documentclass[a4paper,oneside,reqno,11pt]{amsart}
\usepackage{amsmath,amssymb,amsthm}
\theoremstyle{plain}
\newtheorem{theorem}{Theorem}
\theoremstyle{remark}
\newtheorem*{remark}{Remark}

\renewcommand{\d}{\mathrm d}
\newcommand\No{no.~}

\begin{document}

\title[Quadratic transformations and Guillera's formulae]%
{Quadratic transformations\\and Guillera's formulae for $1/\pi^2$}
\author[W.~Zudilin]%
{Wadim Zudilin}
\thanks{The work is supported by a fellowship
of the Max-Planck Institute for Mathematics (Bonn).}
\date{April 2005}
\address{%
\hbox to70mm{\vbox{\hsize=70mm%
\leftline{Department of Mechanics and Mathematics}
\leftline{Moscow Lomonosov State University}
\leftline{Vorobiovy Gory, GSP-2}
\leftline{119992 Moscow, RUSSIA}
\leftline{{\it URL\/}: \texttt{\rm http://wain.mi.ras.ru/}}
}}}
\email{\texttt{wadim@ips.ras.ru}}

\begin{abstract}
We prove two new series of Ramanujan type for $1/\pi^2$.
\end{abstract}

\maketitle

In a series of works \cite{G1}--\cite{G3}, J.~Guillera
discovered several Ramanujan-type formulae for $1/\pi^2$,
some of which he proved. For instance, the following
two identities are proven in~\cite{G1},~\cite{G3}:
\begin{align}
\sum_{n=0}^\infty
\frac{(\frac12)_n^5}{n!^5}
(20n^2+8n+1)\biggl(-\frac1{2^2}\biggr)^n
&=\frac{8}{\pi^2};
\label{eq1}
\\
\sum_{n=0}^\infty
\frac{(\frac12)_n^5}{n!^5}
(820n^2+180n+13)\biggl(-\frac1{2^{10}}\biggr)^n
&=\frac{128}{\pi^2},
\label{eq2}
\end{align}
where $(a)_n=\Gamma(a+n)/\Gamma(a)=\prod_{k=0}^{n-1}(a+k)$
denotes Pochhammer's symbol.
It happens that, in all Guillera's formulae expressions,
the left-hand sides are given by certain hypergeometric
${}_5F_4$-series, with the sole exception~\cite{G4}
\begin{equation}
\sum_{n=0}^\infty A_n\frac{36n^2+12n+1}{2^{10n}}
=\frac{32}{\pi^2},
\qquad\text{where}\quad
A_n=\binom{2n}n^2\sum_{k=0}^n\binom{2k}k^2\binom{2n-2k}{n-k}^2
\label{eq3}
\end{equation}
(this formula remains unproven).
Recall that a generalized hypergeometric series~\cite{Sl}
is defined for $z\in\mathbb C$, $|z|<1$, by
$$
{}_{q+1}F_q\biggl(\begin{matrix}
a_0, \, a_1, \, \dots, \, a_q \\ b_1, \, \dots, \, b_q
\end{matrix}\biggm|z\biggr)
=\sum_{n=0}^\infty
\frac{(a_0)_n(a_1)_n\dotsb(a_q)_n}
{n!(b_1)_n\dotsb(b_q)_n}z^n.
$$

Formula~\eqref{eq3} looks a natural
extension of Sato's formula for $1/\pi$ involving
Ap\'ery's numbers (see \cite{Y1} and \cite{CCL})
as well as of several other examples like
$$
\sum_{n=0}^\infty B_n\frac{4n+1}{36^n}
=\frac{18}{\pi\sqrt{15}},
\qquad\text{where}\quad
B_n=\sum_{k=0}^n\binom nk^4,
$$
proven by Y.~Yang~\cite{Y2}.

The aim of this note is to derive two new identities
of type~\eqref{eq3} from \eqref{eq1},~\eqref{eq2}
using hypergeometric technique.

\section{Quadratic transformations}
\label{s1}

Recall quadratic transformations for ${}_2F_1$-series,
$$
{}_2F_1\biggl(\begin{matrix}
a, \, b \\ 1+a-b
\end{matrix}\biggm|z\biggr)
=(1-z)^{-a}\cdot{}_2F_1\biggl(\begin{matrix}
\frac12a, \, \frac12+\frac12a-b \\ 1+a-b
\end{matrix}\biggm|\frac{-4z}{(1-z)^2}\biggr)
$$
due to Gauss, and for ${}_3F_2$-series,
\begin{equation}
{}_3F_2\biggl(\begin{matrix}
a, \, b, \, c \\ 1+a-b, \, 1+a-c
\end{matrix}\biggm|z\biggr)
=(1-z)^{-a}\cdot{}_3F_2\biggl(\begin{matrix}
\frac12a, \, \frac12+\frac12a, \, 1+a-b-c \\ 1+a-b, \, 1+a-c
\end{matrix}\biggm|\frac{-4z}{(1-z)^2}\biggr)
\label{eq4}
\end{equation}
due to Whipple. Their higher-order analogues necessarily involve
multiple hypergeometric series in the right-hand side
(cf., e.g., \cite{AZ}, Proposition~6, and \cite{Zu}, Theorem~5).
Here we indicate the following result.

\begin{theorem}
\label{th1}
The following quadratic transformation is valid:
\begin{align}
&
{}_5F_4\biggl(\begin{matrix}
a, \, b, \, c, \, d, \, e \\ 1+a-b, \, 1+a-c, \, 1+a-d, \, 1+a-e \end{matrix}
\biggm|z\biggr)
\nonumber\\ &\qquad
=(1-z)^{-a}\sum_{n=0}^\infty\frac{(\frac12a)_n(\frac12+\frac12a)_n}
{(1+a-b)_n(1+a-c)_n}\biggl(\frac{-4z}{(1-z)^2}\biggr)^n
\nonumber\\ &\qquad\quad\times
\sum_{\nu=0}^n
\frac{(b)_\nu(c)_\nu(1+a-d-e)_\nu}
{\nu!(1+a-d)_\nu(1+a-e)_\nu}\,
\frac{(1+a-b-c)_{n-\nu}}{(n-\nu)!},
\label{eq5}
\end{align}
whenever both series converge.
\end{theorem}

\begin{remark}
Theorem~\ref{th1} may be stated in the form of Orr-type theorem
(cf.~\cite{Sl}, Section~2.5):
{\it If
\begin{equation}
(1-z)^{b+c-a-1}
\cdot{}_3F_2\biggl(\begin{matrix}
b, \, c, \, 1+a-d-e \\ 1+a-d, \, 1+a-e
\end{matrix}\biggm|z\biggr)
=\sum_{n=0}^\infty f_nz^n,
\label{eq6}
\end{equation}
then}
\begin{align*}
&
{}_5F_4\biggl(\begin{matrix}
a, \, b, \, c, \, d, \, e \\ 1+a-b, \, 1+a-c, \, 1+a-d, \, 1+a-e \end{matrix}
\biggm|z\biggr)
\\ &\qquad
=(1-z)^{-a}\sum_{n=0}^\infty f_n\frac{(\frac12a)_n(\frac12+\frac12a)_n}
{(1+a-b)_n(1+a-c)_n}\biggl(\frac{-4z}{(1-z)^2}\biggr)^n.
\end{align*}
It follows from \eqref{eq6} that $|f_n|^{1/n}\to1$
as $n\to\infty$, hence the condition $|4z/(1-z)^2|<1$
ensures convergence of the double series in~\eqref{eq5}.
\end{remark}

\begin{proof}
It is sufficient to prove the identity in a neighbourhood of
the origin.
Using the Pfaff--Saalsch\"utz theorem \cite{Sl}, formula (2.3.1.3),
$$
{}_3F_2\biggl(\begin{matrix}
-n, \, a+n, \, 1+a-d-e \\ 1+a-d, \, 1+a-e
\end{matrix}\biggm|1\biggr)
=\frac{(-d-n+1)_n(e)_n}{(1+a-d)_n(e-a-n)_n}
=\frac{(d)_n(e)_n}{(1+a-d)_n(1+a-e)_n},
$$
write
\begin{align*}
&
{}_5F_4\biggl(\begin{matrix}
a, \, b, \, c, \, d, \, e \\ 1+a-b, \, 1+a-c, \, 1+a-d, \, 1+a-e \end{matrix}
\biggm|z\biggr)
\\ &\qquad
=\sum_{n=0}^\infty\frac{(a)_n(b)_n(c)_n}
{n!(1+a-b)_n(1+a-c)_n}z^n\cdot{}_3F_2\biggl(\begin{matrix}
-n, \, a+n, \, 1+a-d-e \\ 1+a-d, \, 1+a-e
\end{matrix}\biggm|1\biggr)
\\ &\qquad
=\sum_{n=0}^\infty\frac{(a)_n(b)_n(c)_n}
{n!(1+a-b)_n(1+a-c)_n}z^n
\sum_{\nu=0}^n\frac{(-n)_\nu(a+n)_\nu(1+a-d-e)_\nu}
{\nu!(1+a-d)_\nu(1+a-e)_\nu}
\\ &\qquad
=\sum_{\nu=0}^\infty\frac{(1+a-d-e)_\nu(-1)^\nu}
{\nu!(1+a-d)_\nu(1+a-e)_\nu}
\sum_{n=\nu}^\infty\frac{(a)_{n+\nu}(b)_n(c)_n}
{(n-\nu)!(1+a-b)_n(1+a-c)_n}z^n
\\ &\qquad
=\sum_{\nu=0}^\infty\frac{(1+a-d-e)_\nu(-1)^\nu}
{\nu!(1+a-d)_\nu(1+a-e)_\nu}
\sum_{m=0}^\infty\frac{(a)_{m+2\nu}(b)_{m+\nu}(c)_{m+\nu}}
{m!(1+a-b)_{m+\nu}(1+a-c)_{m+\nu}}z^{m+\nu}
\displaybreak[2]
\\ &\qquad
=\sum_{\nu=0}^\infty\frac{(a)_{2\nu}(b)_\nu(c)_\nu(1+a-d-e)_\nu}
{\nu!(1+a-b)_\nu(1+a-c)_\nu(1+a-d)_\nu(1+a-e)_\nu}(-z)^\nu
\\ &\qquad\quad\times
{}_3F_2\biggl(\begin{matrix}
a+2\nu, \, b+\nu, \, c+\nu \\ 1+a-b+\nu, \, 1+a-c+\nu
\end{matrix}\biggm|z\biggr).
\end{align*}
Applying the quadratic transformation \eqref{eq4} to
the latter ${}_3F_2$-series, we deduce
\begin{align*}
&
{}_5F_4\biggl(\begin{matrix}
a, \, b, \, c, \, d, \, e \\ 1+a-b, \, 1+a-c, \, 1+a-d, \, 1+a-e \end{matrix}
\biggm|z\biggr)
\\ &\qquad
=\sum_{\nu=0}^\infty\frac{(a)_{2\nu}(b)_\nu(c)_\nu(1+a-d-e)_\nu}
{\nu!(1+a-b)_\nu(1+a-c)_\nu(1+a-d)_\nu(1+a-e)_\nu}(-z)^\nu
\\ &\qquad\quad\times
(1-z)^{-(a+2\nu)}\cdot{}_3F_2\biggl(\begin{matrix}
\frac12a+\nu, \, \frac12+\frac12a+\nu, \, 1+a-b-c \\ 1+a-b+\nu, \, 1+a-c+\nu
\end{matrix}\biggm|\frac{-4z}{(1-z)^2}\biggr)
\\ &\qquad
=(1-z)^{-a}\sum_{\nu=0}^\infty
\frac{(\frac12a)_\nu(\frac12+\frac12a)_\nu(b)_\nu(c)_\nu(1+a-d-e)_\nu}
{\nu!(1+a-b)_\nu(1+a-c)_\nu(1+a-d)_\nu(1+a-e)_\nu}
\biggl(\frac{-4z}{(1-z)^2}\biggr)^\nu
\\ &\qquad\quad\times
\sum_{m=0}^\infty\frac{(\frac12a+\nu)_m(\frac12+\frac12a+\nu)_m(1+a-b-c)_m}
{m!(1+a-b+\nu)_m(1+a-c+\nu)_m}\biggl(\frac{-4z}{(1-z)^2}\biggr)^m
\\ &\qquad
=(1-z)^{-a}\sum_{\nu=0}^\infty
\frac{(b)_\nu(c)_\nu(1+a-d-e)_\nu}
{\nu!(1+a-d)_\nu(1+a-e)_\nu}
\biggl(\frac{-4z}{(1-z)^2}\biggr)^\nu
\\ &\qquad\quad\times
\sum_{n=\nu}^\infty\frac{(\frac12a)_n(\frac12+\frac12a)_n(1+a-b-c)_{n-\nu}}
{(n-\nu)!(1+a-b)_n(1+a-c)_n}\biggl(\frac{-4z}{(1-z)^2}\biggr)^{n-\nu}
\\ &\qquad
=(1-z)^{-a}\sum_{n=0}^\infty\frac{(\frac12a)_n(\frac12+\frac12a)_n}
{(1+a-b)_n(1+a-c)_n}\biggl(\frac{-4z}{(1-z)^2}\biggr)^n
\\ &\qquad\quad\times
\sum_{\nu=0}^n
\frac{(b)_\nu(c)_\nu(1+a-d-e)_\nu}
{\nu!(1+a-d)_\nu(1+a-e)_\nu}\,
\frac{(1+a-b-c)_{n-\nu}}{(n-\nu)!},
\end{align*}
which is the required formula~\eqref{eq5}.
\end{proof}

Plugging in $a=b=c=d=e=\frac12$ we obtain
$$
\sum_{n=0}^\infty\frac{(\frac12)_n^5}{n!^5}z^n
=\frac1{(1-z)^{1/2}}\sum_{n=0}^\infty
\frac{(\frac14)_n(\frac34)_n}{n!^2}
\biggl(\frac{-4z}{(1-z)^2}\biggr)^n
\sum_{\nu=0}^n\frac{(\frac12)_\nu^3}{\nu!^3}\,
\frac{(\frac12)_{n-\nu}}{(n-\nu)!}.
$$
Note the equality
\begin{equation}
u_n=\sum_{\nu=0}^n\frac{(\frac12)_\nu^3}{\nu!^3}\,
\frac{(\frac12)_{n-\nu}}{(n-\nu)!}
=\sum_{\nu=0}^n\biggl(\frac{(\frac14)_\nu(\frac34)_{n-\nu}}
{\nu!(n-\nu)!}\biggr)^2,
\label{eq7}
\end{equation}
followed from \cite{Sl}, formula~(2.5.18).
Another way to prove \eqref{eq7} is based on
the recurrence relation
\begin{equation}
8(n+1)^3u_{n+1}-(2n+1)(8n^2+8n+5)u_n+8n^3u_{n-1}=0
\label{eq8}
\end{equation}
satisfied by both expressions in~\eqref{eq7}.
(Proof of~\eqref{eq8} uses the algorithm
of creative telescoping \cite{PWZ}, Chapter~6.)
We summarize saying above in

\begin{theorem}
\label{th2}
Suppose $|z|<1$ and $|4z/(1-z)^2|<1$.
The following identity holds:
\begin{equation}
\sum_{n=0}^\infty\frac{(\frac12)_n^5}{n!^5}z^n
=\frac1{(1-z)^{1/2}}\sum_{n=0}^\infty
u_n\frac{(\frac14)_n(\frac34)_n}{n!^2}
\biggl(\frac{-4z}{(1-z)^2}\biggr)^n,
\label{eq9}
\end{equation}
where $u_n$ are given in~\eqref{eq7}.
\end{theorem}

\begin{remark}
It is worth mentioning the following curious
hypergeometric identity:
\begin{equation}
\sum_{n=0}^\infty
u_n\frac{(\frac13)_n(\frac23)_n}{n!^2}z^n
={}_3F_2\biggl(\begin{matrix}
\frac16, \, \frac12, \, \frac56 \\ 1, \, 1
\end{matrix}\biggm|z\biggr)^2
={}_2F_1\biggl(\begin{matrix}
\frac1{12}, \, \frac5{12} \\ 1
\end{matrix}\biggm|z\biggr)^4,
\label{eq10}
\end{equation}
since the double hypergeometric series in the left-hand side
is very close to the one used in~\eqref{eq9}.
Formula~\eqref{eq10} was born in our correspondence
with G.~Almkvist and Y.~Yang.
\end{remark}

\section{New formulae for $1/\pi^2$}
\label{s2}

For the differential operator $\theta=z\frac{\d}{\d z}$,
we have
$$
\theta(1-z)^{-1/2}=\frac z{2(1-z)}\cdot(1-z)^{-1/2}
\qquad\text{and}\qquad
\theta\biggl(\frac{-4z}{(1-z)^2}\biggr)
=\frac{1+z}{1-z}\cdot\biggl(\frac{-4z}{(1-z)^2}\biggr).
$$
Therefore,
\begin{align*}
&
\theta\biggl((1-z)^{-1/2}\sum_{n=0}^\infty
C_n\biggl(\frac{-4z}{(1-z)^2}\biggr)^n\biggr)
\\ &\qquad
=(1-z)^{-1/2}\sum_{n=0}^\infty
C_n\cdot\biggl(n\frac{1+z}{1-z}+\frac z{2(1-z)}\biggr)
\cdot\biggl(\frac{-4z}{(1-z)^2}\biggr)^n,
\\ &
\theta^2\biggl((1-z)^{-1/2}\sum_{n=0}^\infty
C_n\biggl(\frac{-4z}{(1-z)^2}\biggr)^n\biggr)
\\ &\qquad
=(1-z)^{-1/2}\sum_{n=0}^\infty
C_n\cdot\biggl(n^2\frac{(1+z)^2}{(1-z)^2}
+n\frac{z(3+z)}{(1-z)^2}+\frac{z(2+z)}{4(1-z)^2}\biggr)
\cdot\biggl(\frac{-4z}{(1-z)^2}\biggr)^n.
\end{align*}
We now apply the results to the functions
$$
20\theta^2f(z)+8\theta f(z)+f(z)
\qquad\text{and}\qquad
820\theta^2f(z)+180\theta f(z)+13f(z),
$$
where $f(z)$ is given in~\eqref{eq9},
and substitute $z=-1/2^2$ and $z=-1/2^{10}$, respectively.
Using then evaluations~\eqref{eq1}, \eqref{eq2}
and formulae
$$
\frac{(\frac14)_n(\frac34)_n}{n!^2}
=2^{-6n}\frac{(4n)!}{n!^2(2n)!},
\qquad
u_n=\sum_{k=0}^n\frac{(\frac12)_k^3}{k!^3}\,
\frac{(\frac12)_{n-k}}{(n-k)!}
=2^{-6n}\sum_{k=0}^n\binom{2k}k^3\binom{2n-2k}{n-k}2^{4(n-k)},
$$
we arrive at

\begin{theorem}
\label{th3}
The following identities are valid:
\begin{align*}
\sum_{n=0}^\infty U_n\frac{(4n)!}{n!^2(2n)!}
\,\frac{18n^2-10n-3}{(2^85^2)^n}
&=\frac{10\sqrt{5}}{\pi^2},
\\
\sum_{n=0}^\infty U_n\frac{(4n)!}{n!^2(2n)!}
\,\frac{1046529n^2+227104n+16032}{(5^441^2)^n}
&=\frac{5^441\sqrt{41}}{\pi^2},
\end{align*}
where the sequence of integers
$$
U_n=\sum_{k=0}^n\binom{2k}k^3\binom{2n-2k}{n-k}2^{4(n-k)},
\qquad n=0,1,2,\dots,
$$
satisfies the recurrence relation
$$
(n+1)^3U_{n+1}-8(2n+1)(8n^2+8n+5)U_n+4096n^3U_{n-1}=0,
\qquad n=1,2,\dotsc.
$$
\end{theorem}

It seems quite likely that the hypergeometric machinery
could lead to several other formulae for $1/\pi^2$.

\medskip
{\bf Acknowledgments.}
I thank Gert Almkvist, Christian Krattenthaler and Yi-Fan Yang
for constructive comments. The final part of the work was
done during my visit of the Max-Planck Institute for Mathematics
in Bonn. I thank the staff of the institute for the wonderful
working conditions I experienced there.


\end{document}